\documentclass[10pt]{article}
\usepackage[applemac]{inputenc}
\usepackage[english]{babel}
\usepackage[T1]{fontenc} 

\usepackage{amsmath} 
\usepackage{amssymb} 
\usepackage{ulem}
\usepackage{amsfonts} 
\usepackage{array} 

\def \N{I\!\!N}
\def \R{I\!\!R}

\def\e{\varepsilon}
\def\di{\displaystyle}

\newtheorem{Theorem}{Theorem}[section]

\newtheorem{Lemma}{Lemma}[section]
\newtheorem{Corollary}{Corollary}[section]

\def\lb{\overline\lambda}
\def\ou{\overline u}
\def\uu{\underline u}
\def\ow{\overline w}
\def\uw{\underline w}
\def\ue{u^{\varepsilon}}

\def\HRL{{\displaystyle{\mathop{\scriptstyle{y\to
x}}_{ t \to +\infty}}}}

\begin{document}
\title{Ergodic type problems and large time behaviour of unbounded solutions of Hamilton-Jacobi Equations}
\author{Guy Barles \\
Laboratoire de Math\'ematiques et Physique \\ Th\'eorique (UMR CNRS 6083) \\ Facult\'e des Sciences et Techniques,  
Universit\'e de Tours \\Parc de Grandmont, 37200 Tours, {\sc France} \\ 
\ \\  
Jean-Michel Roquejoffre, \\
UMR CNRS 5640 and Institut 	Universitaire de France\\Universit\'e Paul Sabatier, 118 route de Narbonne\\31062 Toulouse Cedex, {\sc France}}
\date{ }
\maketitle

\begin{abstract}
We study the large time behavior of Lipschitz continuous, possibly unbounded, viscosity solutions of Hamilton-Jacobi Equations in the whole space $\R^N$. The associated ergodic problem has Lipschitz continuous solutions if the analogue of the ergodic constant is larger than a minimal value $\lambda_{min}$. We obtain various large-time convergence
and Liouville type theorems, some of them being of completely new type. We  also provide  examples showing that, in this unbounded framework, the ergodic behavior may fail, and that the asymptotic behavior may also be unstable with respect to the initial data. 
\end{abstract}

KEY-WORDS : Hamilton-Jacobi Equations, unbounded solutions, ergodic problems, large time behavior, geodesics.


\section{Introduction}  
Recently a lot of works have been devoted to the study of large time behaviour of solutions of Hamilton-Jacobi Equations 
\begin{equation}\label{HJ}
u_t + H(x, Du)=0\quad \hbox{in}\  \R^n \times (0,+\infty),
\end{equation}
\begin{equation}
\label{idHJ}
u(x,0)=u_0(x)\quad \hbox{in}\  \R^n .
\end{equation}
The usual assumptions are: $H \in C (\R^n \times \R^n)$ and $u_0 \in W^{1,\infty}(\R^n)$ are periodic in $x$, while $H(x,p)$ is convex and coercive
in p, i.e.
\begin{equation}
\label{coer}H(x,p)\to +\infty\quad \hbox{as}\ |p|\to +\infty\quad \hbox{uniformly w.r.t.}\ x \in \R^n.
\end{equation}
As a consequence of these assumptions, the solutions of (\ref{HJ})-(\ref{idHJ}) are Lipschitz continuous and periodic in $x$ and, in good cases, they are expected to remain uniformly bounded in $x$ and $t$ and to converge uniformly as $t\to +\infty$ to solutions of the stationary equation which are also Lipschitz continuous and periodic in $x$. In particular, a key feature in these results is the boundedness of the solutions, both of the evolution and stationary equations. A notable exception to this are the papers of Contreras \cite{C} and Fathi \& Maderna \cite{FaMad}, where the periodicity assumption is dropped, and where the existence of possibly unbounded solutions of the stationary equation is looked for. 

The aim of the present paper is to present a systematic study of cases where one has non-periodic and - this is the main point - unbounded solutions, in particular for the limiting stationary equation. In the periodic setting, one first solves a so-called ergodic problem, namely a stationary
Hamilton-Jacobi Equation of the type
\begin{equation}
\label{EP}H(x,D\ou)=\lambda \quad \hbox{in}\ \R^n.
\end{equation}
where both the function $\ou$ and the constant
$\lambda$ are unknown.
From Lions, Papanicolaou \& Varadhan\cite{LPV}, there exists a unique constant $\lambda = \lb$ such that (\ref{EP}) has a Lipschitz continuous, periodic solution. It is worth remarking that the actual interest of this result is to
produce a bounded solution, and this is where periodicity
plays a key role.
The connection with large time
behaviour in (\ref{HJ}) is
then the following : on the one hand, one can prove that the
solution $u$ of (\ref{HJ})-(\ref{idHJ}) satisfies
\begin{equation}\label{ust}\frac{u(x,t)}t \to -\lb \quad \hbox{as}\ t \to +\infty\quad \hbox{uniformly in}\ \R^n \end{equation}
and, on the other hand, that
\begin{equation}
\label{complt}u(x,t)+ \lb t \to u_\infty(x)\quad \hbox{as $t$} \to +\infty\quad \hbox{uniformly in}\ \R^n\,,
\end{equation}
where $u_\infty$ is a solution of (\ref{EP}) with $\lambda = \lb$.
It is worth pointing out that, if a property like (\ref{ust}) can be obtained rather easily as a consequence of standard comparison results for equation (\ref{HJ}), the more precise asymptotic behaviour (\ref{complt}) is, on the contrary, a far more difficult result~; in fact, the asymptotic behaviour of solutions of (\ref{HJ})-(\ref{idHJ}) remained an open problem for a long
 time.
 Namah \& Roquejoffre \cite{NR} were the first to break this difficulty
 under the following additional assumptions
 \begin{equation}\label{hypNR}H(x,p) \geq H(x,0)\quad \hbox{in}\ \R^n \times \R^n \quad\hbox{and} \quad\ \max_{\R^n} H(x,0)=0.
 \end{equation}
 This assumption seems to be a bit restrictive
 but, on one hand, it covers several interesting cases and, on the
 other hand, this result does not require strong convexity
 assumptions on $H$ in $p$.
 
Then a second type of result was obtained by Fathi \cite{Fa} whose proof was based on dynamical systems type arguments and in particular on the so-called Mather set which is (roughly
speaking) an attractor for the geodesics associated to the representation formula of $u$. Contrarily to \cite{NR}, the results of \cite{Fa} use rather strong convexity assumptions on $H$ (and
also far more regularity) but do not need (\ref{hypNR}). In fact, (\ref{hypNR}) can be interpreted,
in the strictly convex case, as a  special assumption on the Mather set. Fathi's results were extended to time-dependent hamiltonians in one space variable in \cite{BR}.

The most general result is this direction is the one of Barles \& Souganidis \cite{BS1}, which generalizes both results of \cite{NR} and \cite{Fa}, and which can even handle some special cases where $H$ is not convex in $p$. The key assumption is on the quantity $H_p(x,p)\cdot p - H(x,p)$ which, in some sense, measures the attractivity of the Mather set. Natural questions are then~: can such results be obtained without assuming periodicity? Are they true (with some natural modifications of statements) for unbounded Lipschitz continuous solutions? Are these results "stable"
under (non-periodic) small pertubations? A negative answer is given by a counter-example due to Barles \& Souganidis \cite{BS2}, which shows that the above results are wrong if one drops the periodicity assumption on $u_0$, even if $u_0$ remains bounded and Lipschitz continuous. This counter-example emphasizes that the behaviour at infinity of $u_0$ may play a role to deduce the behaviour of $u$ as $t \to  +\infty$.

It is to be noted that such ergodic problems arise also in homogenization theory (the so-called ``cell problems'') and a related question to ours is whether the periodicity assumption on $H$ can be removed while keeping bounded solutions (the ``correctors''): we refer to Ishii\cite{I2}  for the existence of bounded approximate correctors in the almost periodic framework and to Lions \& Souganidis \cite{LS} for a complete discussion of this problem, not only in the deterministic framework but also for equations with a stochastic dependence. 

The main results of our paper is that convergence results survive under more stringent assumptions and, if we insist on weakening the assumptions on $u_0$ as much as possible, Liouville type theorems are still available. To summarize, we prowe the following. 
\begin{itemize}
\item[(i)] Under assumption (\ref{coer}) and if $H$ is bounded, uniformly continuous on $\R^{n}\times B(0,R)$ for any $R>0$, there exists $\lambda_{min}\in\R$ such that the ergodic problem \eqref{EP} has solutions if and only if $\lambda\geq\lambda_{min}$.
\item[(ii)] If we assume, in addition, that $H(x,p)$ is convex in $p$ and that $u_0(x)-\phi(x) \to 0$ as $|x| \to + \infty$, where $\phi$ is a solution of the ergodic problem for some $\lambda>\lambda_{min}$, then the solution $u$ of the Cauchy Problem \eqref{HJ}-\eqref{idHJ} satisfies
$$
\lim_{t\to+\infty}(u(x,t)-\lambda t-\phi(x))=0\ \ \ \hbox{locally uniformly in $\R^n$}.
$$
\item[(iii)] Under suitable additional assumptions on $H$ of strong convexity type, if $u$ is a solution of \eqref{HJ} in $\R^n\times\R$ such that there exists a sub-solution $\underline\phi$ of the ergodic problem for some $\lambda\geq\lambda_{min}$ for which $u(x,t)-\lambda t-\underline\phi(x)$ remains bounded on $\R^n\times\R$, then there is a solution $\ou$ of \eqref{EP} such that $u(x,t)=\lambda t+\ou(x)$.
\end{itemize}

We complement these positive results by describing various pathologies arising when the boundedness assumptions on the solution of \eqref{HJ} is removed. In particular, even the ergodic behavior may fail as is shown in Section 3. 

\medskip
The present paper is organized as follows : in Section 2 we state
 and prove the result concerning the solutions of (\ref{EP}).
 Section 3 is devoted to the description of various troubles encountered
 in the unbounded context: loss of stability and uniform convergence, loss of the property (\ref{ust}).
In Section 4, we provide the results on the convergence of the solution of (\ref{HJ})-(\ref{idHJ}) as $t \to +\infty$, thus covering Point (ii) above. Finally, we prove the Liouville-type result - Point (iii) above - in Section 5.

\section{Bounded and unbounded solutions of the ergodic
 equations}

The following theorem is proved in Fathi \& Mather\cite{FaMa} in the strictly convex case, using the Lax-Oleinik formula. We provide here an alternative proof, also valid in the nonconvex case.
\begin{Theorem}
 \label{ErgRes}
Assume that $H$ is bounded, uniformly continuous on $\R^{n}\times B(0,R)$ for any $R>0$ and that (\ref{coer}) holds. Then there exists $\lambda_{min} \in \R$ such that, for any $\lambda \geq \lambda_{min}$, there exists a Lipschitz continuous solution of (\ref{EP}).
  \end{Theorem}
 
 \noindent{\bf Proof.}
{\bf 1.} We first prove that, if 
 $$
 \lambda>\sup_{x\in\R^n}H(x,0),
 $$
then (\ref{EP}) has a Lipschitz continuous solution. To see this, we first notice that $0$ is a subsolution of the equation. Then we consider $R>0$ and the Dirichlet problem
 \begin{equation}
 \label{EPR}
 H(x,Du)=\lambda\ \ \hbox{in $B_R(0)$,}\ \ \ u=0\ \hbox{on $\partial B_R(0)$.}
 \end{equation}
If $C_R>0$ is large enough and the vector $p\in\R^n$ has a large enough norm, then the function $x\mapsto C_R+p.x$ is a positive super-solution to (\ref{EPR}). Consequently, by the Perron's method, combining classical arguments of Ishii\cite{I} (see also \cite{GBlivre}) and the version up to the boundary of Da Lio\cite{FDL}, one easily shows that (\ref{EPR}) has a Lipschitz continuous solution that we call $u_R$. Then the function
$$
  v_R=u_R-u_R(0)
  $$
vanishes at 0 and, by (\ref{coer}), its gradient is uniformly bounded in $R$. Using Ascoli's theorem together with the classical stability result for viscosity solutions gives the convergence of a subsequence $(v_{R_n})_n$ to a solution of (\ref{EP}).

\noindent{\bf 2.} Denote by $\lambda_{min}$ the infimum of all $\lambda$ such that (\ref{EP}) has solutions. We claim that $\lambda_{min}$ is not $-\infty$: indeed, any solution of (\ref{EP}) satisfies, almost everywhere:
  $$
  H(x,Du)\geq\inf_{(x,p)\in\R^n\times\R^n}H(x,p).
  $$
Consequently, $\lambda_{min}$ has to be larger than the above right-hand side.

\noindent{\bf 3.} Let us prove that (\ref{EP}) has a solution for $\lambda=\lambda_{min}$. Without loss of generality we may assume the existence of a sequence $(\lambda_n)_n$ converging to $\lambda_{min}$ for which there is a solution $u_n$ to (\ref{EP}). Then the family $(v_n)_n$ given by $v_n=u_n-u_n(0)$ is relatively compact in $C(\R^n)$ and using again Ascoli's theorem together with the classical stability result for viscosity solutions yields a solution $u_{\lambda_{min}}$ to (\ref{EP}) for $\lambda=\lambda_{min}$.

\noindent{\bf 4.} In order to conclude that (\ref{EP}) has solutions for all $\lambda\geq\lambda_{min}$, we repeat exactly the argument of Step~1 above, except for a slight point~: instead of using $0$ as a subsolution, we use $u_{\lambda_{min}}$ and we replace in (\ref{EPR}), the boundary condition by ``$u=u_{min}$ on $\partial B_R(0)$''.\hfill $\bullet$

\medskip
At that point, it is worth making the following comment: if we assume that $H$ is periodic in $x$, then, as pointed out at the beginning of the introduction, we know from \cite{LPV} that there is $\lb$ such that the ergodic problem has a bounded and periodic solution if and only if $\lambda=\lb$. We notice here that there is no reason why we should have $\lb=\lambda_{min}$; indeed we always have $\lb \geq \lambda_{min}$, but the strict inequality may hold: indeed, consider  in one space dimension
$$
H(x,p)=\vert p - 1\vert.
$$
Then we have $\lambda_{min}=0$ - just because $x\mapsto x$ solves the ergodic problem with $\lambda=0$ and clearly $\lambda_{min} \geq 0$- and 
$\lb=1$ - simply because $x\mapsto 0$ is periodic in $x$ and solves the ergodic problem with 
$\lambda=0$. We refer to \cite{FaMad} for a related study.

\section{Some pathologies of the unbounded setting}
We analyze in this section various troubles occurring in the non-periodic, unbounded setting. 
A first example - constructed on the Barles \& Souganidis model \cite{BS1} shows that, in the unbounded setting, ergodic behaviour is very easily lost. In a second paragraph we study some instabilities with respect to the 
Hamiltonian. Such instabilities are already present in the periodic setting, but the very strong convergence property makes them less visible. 
\subsection{A counter-example to the ergodic behaviour}
  The counter-example is provided in the following
  \begin{Theorem} 
  There exists a Lipschitz continuous
  initial data $u_{0}$ in $\R$, such that, if $u$ is the solution of
  \begin{equation}
  \label{cexHJ}
  u_t -u_{x} + \frac1{2}|u_{x}|^2=0\quad \hbox{in}\ \R \times (0,+\infty)\; ,
  \end{equation}
  then 
  $\displaystyle\frac{u(0,t)}{t}$ 
  does not converge as $t\to +\infty$.
  \end{Theorem}

  \noindent{\bf Proof.} The solution of (\ref{cexHJ}) associated with a Lipschitz continuous initial data $u_0$ is given by the Oleinik-Lax formula
  $$u (x,t)=\displaystyle\inf_{y \in \R}\biggl(u_{0}(y) + \frac{|x+t-y|^2}{2t}\biggr)\; .
  $$
  Let $(a_n)_{n\in \N}$ be a strictly increasing sequence of non-negative real numbers such that
  \begin{equation}\label{propsuit}
  \lim_n\,\frac{a_{n+1}}{a_n} = +\infty\; .
 \end{equation}
   We consider the Lipschitz continuous initial data $u_0$ defined in the following way
   $$u_0(y) = 0 \quad\hbox{for  }y\leq a_0\; ,$$
   and for any $k\in \N$
   $$ u'_0 (y) = 
   \left\{\begin{array}{ll}
   0 & \hbox{if  }y \in (a_{2k+1}, a_{2k+2})\; ,\\
   -1 & \hbox{if  }y \in (a_{2k+2}, a_{2k+3})\; .
   \end{array}\right.
   $$
   Now we examine $u(0,t)$. Since $-1 \leq u'_0(y) \leq 0$ in $\R$, one checks easily that the infimum in the Oleinik-Lax formula is achieved at $\overline y$ which satisfies
   $$ t \leq \overline y \leq 2t\; .$$
   For $k\in \N$ large enough, we first consider the case when 
   $t \in \displaystyle (a_{2k+1}, \frac{1}{2}a_{2k+2})$~: 
   since $u_0$ is constant on this interval and taking account of the property of 
   $\overline y$ above, one has clearly $\overline y=t$ and therefore
   $$ u(0,t) = u_0(a_{2k+1}).
   $$
   Using this for $\displaystyle t_k= \frac{1}{4}a_{2k+2} > a_{2k+1}$ (we recall that (\ref{propsuit}) holds and that $k$ is chosen large enough), we deduce
   $$\frac{u(0,t_k)}{t_k} = \frac{4u_0(a_{2k+1})}{a_{2k+2}} \to 0\quad\hbox{as  }k\to \infty\; .
   $$
   Indeed, since $u_0$ is Lipschitz continuous with a Lipschitz constant equal to $1$, 
   $|u_0(a_{2k+1})|\leq a_{2k+1}$ and the above property is a consequence of the choice of the sequence
    $(a_n)_{n\in \N}$.
    Now we perform the same argument but for $t$ in intervals of the form 
    $ \displaystyle (a_{2k}, \frac{1}{2}a_{2k+1})$. This time, the optimization provides
    $$ \overline y = 2t\; ,$$
    and
    $$ u(0,t) = u_0 (2t) + \frac{t}{2}\; .$$But $u_0 (2t) = u_0(a_{2k}) - (2t-a_{2k})$ 
    and therefore by choosing $t'_k= \frac{1}{4}a_{2k+1}$ (again $t'_k > a_{2k}$ by (\ref{propsuit}) and the fact that $k$ is chosen large enough), we have
    $$ \frac{u(0,t'_k)}{t'_k} = \frac{1}{t'_k}\left(u_0(a_{2k}) - (2t'_k-a_{2k})+ \frac{t'_k}{2}\right) \to -\frac{3}{2}\; ,$$
    by using again the main properties of $u_0$ and the sequence 
    $(a_n)_{n\in \N}$.
    Therefore we have two different limits for the sequences
     $\displaystyle\left(\frac{u(0,t_k)}{t_k}\right)_k$ and $\displaystyle\left(\frac{u(0,t'_k)}{t'_k}\right)_k$ 
     with $t_k, t'_k \to +\infty$, and the counter-example is complete.\hfill $\bullet$
     
\subsection{Instability with respect to the initial data}

Let us formulate the following very simple question~: under  ``good conditions'' on $H$ and $u_0$, what can we say about the large time behaviour of the
solution $\ue$ of
 $$\ue_t + H(x,D\ue)=\varepsilon f(x)\quad \hbox{in}\ \R^n \times (0,+\infty)\, ,$$ 
 $$\ue(x,0)=u_0(x)+ \varepsilon g(x)\quad \hbox{in}\ \R^n\, ,
 $$
 where, say, $f,g$ are $C^\infty$-function with compact supports
 and $\varepsilon \ll 1$? Is there some stability with respect to the initial data and the right-hand side of the equation?
 
 We denote by $\varphi$ a $C^\infty$ function with compact
 support such that $\displaystyle\min_{\R^n} \varphi =\varphi(0)=-1$
 and we first consider the case $f = 0,g= \varphi$ and $u_0 \equiv0$.
 If we consider the Hamilton-Jacobi Equation
 $$\ue_t + \frac{1}{2}|D\ue|^2=0\quad \hbox{in}\ \R^n \times (0,+\infty)$$
 then, by the Oleinik-Lax formula, $\ue$ is given by 
 $$\ue(x,t)=\displaystyle\inf_{y\in \R^n}\biggl(\varepsilon g(y)+\frac{|x-y|^2}{2t}\biggr)
 $$
 and it is easy to see that $\ue(x,t) \to - \varepsilon$ locally uniformly while, for any $t, \ue(x,t) \to 0 $ as $|x| \to +\infty$. In this case, the perturbation has a (slight)
 effect and changes a little bit the asymptotic behaviour of the solution.

If, on the other hand, we consider the pde
 $$\ue_t -e\cdot D\ue + \frac{1}{2}|D\ue|^2=0\quad \hbox{in}\ \R^n \times (0,+\infty)
 $$ 
where $e\in \R^n-\{0\}$,  then the solution is given by
 $$\ue(x,t)=\displaystyle\inf_{y \in \R^n}\biggl(\varepsilon g(y)+ \frac{|x+te-y|^2}{2t}\biggr)
 $$
 and, this time, $\ue(x,t)\to 0$ locally
 uniformly as $t \to +\infty$, while $\ue(-te,t) \equiv - \varepsilon$. 
 Here the behaviour seems to be the same
 as it was without the perturbation but we loose anyway
 again the uniform convergence in $\R^{n}$.

 These two examples show that the effects of the perturbation
 can be rather different (depending on $H$) but, in both
 cases, the uniform convergence in $\R^n$ as $t \to+\infty$ cannot be true anymore and 
 one has to switch to
 a local convergence type requirement. 
 Unfortunately we are
 unable to provide any general result in this direction.
 Moreover we can point out that if, in the
 second example above, we remove the assumption that $g$
 has a compact support then we are exactly in the setting of the counter-example of Barles \& Souganidis \cite{BS2} and therefore
 we do not have convergence anymore.
 The effect of the perturbation $\varepsilon f$ is even stronger :
 to show this, let us consider now the case when $f =\varphi$, $g \equiv 0$ and the pde is the following
 \begin{equation}\label{pert}\ue_t + |D\ue|^2 = \varepsilon f(x) \quad \hbox{in} \ \R^n \times (0,+\infty).\end{equation}
 The problem is here that, if we consider the stationary equation
 $$|D\ou|^2 = \varepsilon f(x)+ \lambda \quad \hbox{in} \ \R^n,$$
 then there is no $\lambda$ for which this pde has a bounded
 solution. This is a striking difference with the Lions, Papanicolaou \& Varadhan
 result and this shows that there is no hope to have a result
 like (\ref{complt}) with a bounded $u_\infty$. Fortunately, here, if we choose
 $\lambda=\varepsilon$, an approach of the type \cite{NR}
 applies and we are able to show that $\ue(x,t) - \varepsilon t\to u_\infty (x)$ as $t \to +\infty$, locally uniformly where $u_\infty$ is a solution of
 $$|Du_\infty|^2=\varepsilon f(x) + \varepsilon\quad \hbox{in}\ \R^n.$$
This result is a consequence of Theorem~\ref{t4.1} in Section~\ref{AR}.

\subsection{Outline of the rest of the paper}

 We examine in the rest of the paper the large time behaviour of solutions of (\ref{HJ})-(\ref{idHJ}) i.e. the validity of a property like (\ref{complt}); again we consider the case when $u_0$ is a Lipschitz continuous, possibly unbounded, function and of course "uniformly" has again to be replaced by "locally uniformly". We obtain in this direction two types of results for convex Hamiltonians which are in some sense complementary~: the first one is a generalization of the result of \cite{NR} in this non-periodic and even unbounded framework : here we need $u_0$ to be bounded from below for reasons explained below. As a consequence we can analyse completely equation (\ref{pert}).

The second type of result is more original : we assume that $H$ is convex and $u_0(x)-\phi(x) \to 0$ at infinity where $\phi$ is solution of (\ref{EP}) for some $\lambda$.
We prove that, if $\lambda > \lambda_{min}$, then $u(x,t)+\lambda t\to \phi (x)$ locally uniformly as $t \to +\infty$. Therefore, in this case, the large time behaviour of solution is governed by the behaviour for large $x$ of the initial data: we point out that both the ``$\lambda$'' which is selected and the limit of $u(x,t)+\lambda t$ depends on
$\phi$, i.e. on the behaviour of $u_{0}$ for large $x$. Such a behaviour was already observed but with a far less generality in Barles \cite{gb1}. In the case when $H$ satisfies an assumption of the type \eqref{hypNR}, this behaviour shows, on the one hand, that $\lambda_{min}$ is the only constant for which (\ref{EP}) has a solution which is bounded from below, and, on the other hand, it justifies the assumption ``$u_0$ bounded from below'' made in Theorem~\ref{t4.1} below: indeed in this case, the behaviour is always governed by $\lambda_{min}=0$.

The interpretation of this result is rather clear from its proof~: for $\lambda >\lambda_{min}$,  the geodesics have to go to infinity. This is why in this framework, the behaviour of $u_0$ at infinity plays a key role in the determination of the behaviour of  $u$ as $t \to +\infty.$ This is completely different under Condition \eqref{hypNR}, where the geodesics are attracted by the compact set $K:=\{x\in \R^n;\ H(x,0)=0\}$.

 \section{Large-time convergence}\label{AR}

 \begin{Theorem} {\bf (Unbounded version of Namah \& Roquejoffre \cite{NR})}
  \label{t4.1}
Under the assumptions of Theorem~\ref{ErgRes}, if $u_0$ is a bounded from below, Lipschitz continuous function and if $H$ is convex in $p$ and satisfies
      \begin{equation}
      \label{hypNRbis}
      H(x,p) \geq H(x,0)\quad \hbox{in}\ \R^n \times \R^n \;,
      \end{equation}
with $\displaystyle \max_{\R^n} H(x,0)=0$, the set $K:=\{x\in \R^{n}; H(x,0)=0\}$ is a non-empty compact subset of $\R^{n}$ and
\begin{equation}\label{condinf}
\limsup_{|x|\to +\infty}H(x,0)<0\; ,
\end{equation}
then the solution $u$ of (\ref{HJ})--(\ref{idHJ}) converges as $t\to +\infty$
to a solution of (\ref{EP}) with $\lambda = \lambda_{min}=0$.
\end{Theorem}

Before providing the proof of this result, we complement it by the

\begin{Theorem}
      \label{t4.2}
Assume that the assumptions of Theorem~\ref{ErgRes} hold and that $H$ is convex in $p$. If the initial data $u_0$ satisfies
\begin{equation}
\label{idinf}
\lim_{|x| \to +\infty}\, (u_0(x) -\phi(x)) = 0\; ,
\end{equation}
where $\phi: \R^n \to \R$ is a solution of (\ref{EP}) for some $\lambda >\lambda_{min}$, then we have
$$ u(x,t) + \lambda t \to \phi(x) \quad \hbox{locally uniformly in $\R^n$ as $t\to +\infty$.}$$
\end{Theorem}

\noindent{\bf Proof of Theorem~\ref{t4.1}.} 
{\bf 1.} We start by some basic estimates. Since $u_0$ is bounded from below, we can consider  $M=||(u_0)^-||_\infty$ and since $u_0$ is Lipschitz continuous we can introduce its Lipschitz constant $K$. We notice that $-M$ is a subsolution of (\ref{HJ}), while for $x_0 \in K$ and C large enough, $C|x-x_0| + C$ is a supersolution of (\ref{HJ}). By choosing in particular $C>K$, we have
$$ -M \leq u_0(x) \leq C|x-x_0| + C \quad \hbox{in  }\R^n\; ,$$
and, by the maximum principle, we have
$$ -M \leq u(x,t) \leq C|x-x_0| + C \quad \hbox{in  }\R^n\times (0,+\infty)\; .$$
On the other hand, we also have - see \cite{NR} for a proof:
$$ |u_t (x,t)|, |Du(x,t)| \leq \tilde C \quad \hbox{in  }\R^n\times (0,+\infty)\; ,$$
for some large enough constant $\tilde C$ depending only on $H$ and $u_0$.

\noindent{\bf 2.} Using similar sub and supersolutions and repeating the argument of the proof of Theorem~\ref{ErgRes}, we see that one has a solution of (\ref{EP}) for $\lambda = 0$ and therefore $\lambda_{min}\leq 0$. But, for $\lambda < 0$, no solution can exists since $H(x,p)-\lambda >0$ on $K$, therefore $\lambda_{min} = 0$.

\noindent{\bf 3.} On the compact set $K$, $H(x,p) \geq 0$ for any $p$ and therefore $u(x,t)$ is a decreasing function of $t$. This implies the uniform convergence of $u$ to a continuous function $\varphi$; we refer to \cite{NR} for a more detailed proof of this fact.

\noindent{\bf 4.} On $\R^n\backslash K$, we use the half-relaxed limit method and introduce
$$ \ou (x) : = \limsup_{\HRL}\,u (y)\quad , \quad \uu (x) : = \liminf_{\HRL}\,u (y)\; .$$
These functions are respectively sub and supersolutions of the Dirichlet problem
$$ H(x,Dw) = 0\quad\hbox{in  }\R^n-K\; ,$$
$$ w = \varphi \quad\hbox{on  }K\; .$$
It is worth pointing out that, because of the estimates of Step 1, $\ou$ and $\uu$ are Lipschitz continuous on $\R^n$; we also have $\ou\geq\uu$ and $\ou=\uu=\varphi$ on $K$. 

\noindent{\bf 5.} The final point consists in comparing $\ou$ and $\uu$. The fact that the constants are strict sub-solutions does not seem to apply easily here due to the unboundedness of the domain. We use instead a remark of Barles \cite{gb2} (See also \cite{GBlivre}, p. 40): for a given closed bounded convex set $C$ with nonempty interior and containing 0 in its interior, consider its gauge - with respect to 0 - $j_C(p)$ defined as
$$
j_C(p)=\inf\{\lambda>0:\ \frac{p}\lambda\in C\}.
$$
We have $p\in C$ if and only if $j_C(p)\leq 1$, and $p\in\partial C$ if and only if $j_C(p)=1$.

For $\varepsilon>0$ small, we are going to argue in the domain $O_\varepsilon:=\{x:\; H(x,0)<-\varepsilon\} $. Because of condition (\ref{condinf}), if $\varepsilon$ is small enough, the $\partial O_\varepsilon$ remains in a compact subset of $\R^n$ and, for any $x \in \partial O_\varepsilon$, $d(x,K)\leq \rho(\varepsilon)$ where $\rho(\varepsilon)\to 0$ as $\varepsilon \to 0$.

In $O_\varepsilon$, since $0$ is in the interior of the convex set
$$
C(x)=\{p\in\R^n:\ H(x,p)\leq 0\},
$$ 
we can transform the equation $H(x,Dw)=0$ into $G(x,Dw) = 1$, where
$$
G(x,p)=j_{C(x)}(p).
$$ 
The function $G$ satisfies the same assumptions as $H$ and is is also homogeneous of degree $1$ in $p$. Then we may use the Kruzhkov's transform
$$ \ow (x):= -\exp(-\ou(x))\quad , \quad  \uw (x):= -\exp(-\uu(x))\; .$$
The functions $\ow$ and $\uw$ are respectively sub and supersolutions of 
$$ G(x,Dw) + w = 0\quad\hbox{in  }\R^n-K\; .$$
Moreover, $\ow$ and $\uw$ are bounded and even Lipschitz continuous. 

Finally, on $\partial O_\varepsilon$, we have $\ou - o_\varepsilon (1)\leq \varphi \leq \uu+ o_\varepsilon(1)$ by the above mentioned property on $\partial O_\varepsilon$ and this yields $\ow \leq \uw + o_\varepsilon(1)$ on $\partial O_\varepsilon$.

A standard comparison result then applies - see \cite{CIL}, \cite{LiLivre} - and shows that $\ow \leq \uw+ o_\varepsilon(1)$ in $O_\varepsilon$. By letting $\varepsilon$ tends to zero, we obtain that $\ow \leq \uw$ in $\R^n-K$ and therefore the same inequality holds for $\ou$ and $\uu$. By standard arguments, this implies the local uniform convergence of $u$ to the continuous function $u_\infty:=\ou=\uu$ in $\R^n$.\hfill$\bullet$ 

\medskip

\noindent{\bf Proof of Theorem~\ref{t4.2}.} We prove this result in the case when $H$ is superlinear in $p$ and when $L$, the Lagrangian associated to $H$, is also superlinear in $p$ since this case contains most of the interesting ideas. The other cases follow from suitable (easy) adaptations of the arguments, in particular by changing the type of Oleinik-Lax formula we are going to use below.

We recall that $L$ is given, for $x\in \R^n$ and $v\in \R^n$ by
\begin{equation}
      \label{e4.1}
      L(x,v)=\inf_{p\in\R^n}(p.v-H(x,p))\; ,
      \end{equation}
and that the solution $u$ is given by the Oleinik-Lax formula
\begin{equation}
      \label{e4.2}
      u(t,x)=\inf_{\gamma(t)=x}\biggl(u_0(\gamma(0))+\int_0^tL(\gamma(s),\dot\gamma(s))\ ds\biggl),
      \end{equation}
the infimum being taken on the space of absolutely continuous paths $\gamma$ such that $\gamma(t)=x$. We point out that the first simplification in the additional assumptions we made above is that this formula takes such a simple form since, in particular, $L$ is finite for any $x$ and $v$. 

This infimum (and this is where the superlinearity of $L$ plays a role) is attained for an absolutely continuous curve $(\gamma_t (s))_{s\in[0,t]}$\footnote{In more general cases, one may just use approximate minimizers.}.

The proof of Theorem \ref{t4.2} relies on a lemma which is almost as important as the theorem itself.
 \begin{Lemma}
     \label{l4.1}
     Under the assumptions of Theorem \ref{t4.2}, for any $x\in \R^N$, we have
     \begin{equation}
     \label{e4.3}
     \lim_{t \to+\infty}\vert\gamma_{t}(0)\vert=+\infty.
     \end{equation}
\end{Lemma}
     
Let us notice that this lemma implies the following statement of independent interest: if $\lambda>\lambda_{min}$, then there is no bounded extremals associated to a solution $\phi$ of (\ref{EP}), even though there might be bounded solutions - for instance in the periodic setting. This is a striking difference with the Namah-Roquejoffre case where the set $K$ attracts the geodesics.
     
 \medskip
\noindent{\bf Proof.} We first remark that $-\lambda t + \phi$ is a solution of the evolution equation; therefore by the contraction principle
\begin{equation}\label{est1}
||u(x,t)+\lambda t -\phi(x)||_\infty \leq ||u_0-\phi||_\infty\; ,
\end{equation}
and since the right-hand side of (\ref{est1}) is finite by (\ref{idinf}), we deduce that the function $u(x,t)+\lambda t -\phi(x)$ is uniformly bounded.

We may assume, without loss of generality, that $\lambda=0$ and $\lambda_{min}<0$. Because of Theorem~\ref{ErgRes}, for every small enough $\e>0$, there is a solution to (\ref{EPR}) with $\lambda=-\e$. We choose such a $\e$ and denote by $\phi_{-\e}$ a corresponding solution.

We assume, by contradiction, that the lemma is false and that there exists a sequence $(t_n)_n$ converging to $+\infty$ and such that $\gamma_{t_n}(0)$ remains bounded. Since $\e t + \phi_{-\e}$ is a solution of (\ref{HJ}), by the Oleinik-Lax formula,  we have
$$ \e t_n + \phi_{-\e}(x)=\inf_{\gamma(t_n)=x}\biggl(\phi_{-\e} (\gamma(0))+\int_0^{t_n} L(\gamma,\dot\gamma)ds\biggl),$$
while, by the optimality of $\gamma_{t_n}$
$$
 u(x,t_n)=u_0(\gamma_{t_n}(0))+\int_0^{t_n}L(\gamma_{t_n},\dot\gamma_{t_n})ds\; .$$
Therefore
$$
\begin{array}{rll}
     \e t_n+\phi_{-\e}(x)\leq&\phi_{-\e}(\gamma_{t_n}(0))+\di\int_0^{t_n}L(\gamma_{t_n},\dot\gamma_{t_n})\ ds\\
     =&\phi_{-\e}(\gamma_{t_n}(0))-u_0(\gamma_{t_n}(0))+u(t_n,x)
     \end{array}
$$
This property is a contradiction for $n$ large enough since the left-hand side tends to infinity with $n$, while the right-hand side remains bounded because of the assumption on $\gamma_{t_n}(0)$ for the two first terms and the estimate (\ref{est1}) for the last one. \hfill$\bullet$
     
\medskip
 We come back to the proof of Theorem \ref{t4.2}. For $\e>0$, by (\ref{idinf}), there exists $\rho_\e>0$ such that
     \begin{equation}
     \label{e4.6}
     \sup_{|x|\geq \rho_\e}\vert u_0(x)-\phi(x)\vert\leq\e.
     \end{equation}
On the other hand, from Lemma \ref{l4.1}, there is $t_\e>0$ such that, for $t\geq t_\e$, Formula (\ref{e4.2}) becomes
\begin{equation}
     \label{e4.4}
u(t,x)=\di\inf_{\vert\gamma(0)\vert\geq \rho_\e }\biggl(u_0(\gamma(0))+\int_0^tL(\gamma,\dot\gamma)\ ds\biggl).
     \end{equation}
Similarly, by applying Lemma~ \ref{l4.1} to the solution $\phi - \lambda t$, we have
    \begin{equation}
     \label{e4.5}
     -\lambda t+\phi(x)= \inf_{\vert\gamma(0)\vert\geq \rho_\e }\biggl(\phi(\gamma(0))+\int_0^tL(\gamma,\dot\gamma)\ ds\biggl).
     \end{equation}
 Combining (\ref{e4.6}) and (\ref{e4.4}) together with the property $\vert\inf(\cdots)-\inf(\cdots)\vert \leq \sup \vert \cdots - \cdots\vert$, yields
$$
  \vert u(x,t)+ \lambda t-\phi(x)\vert 
 \leq  \sup_{\vert \gamma(0)\vert \geq \rho_\e}  \vert \phi(\gamma(0))-u_0(\gamma(0))\vert \leq 2\e\,.
$$
This provides the pointwise convergence. But since, for $t>0$, the function $x\mapsto u(x,t)+\lambda t -\phi(x)$ is in a compact subset of $C(\R^n)$, this pointwise convergence implies the local uniform convergence. \hfill $\bullet$ 

We notice the following consequence of Theorem \ref{t4.2}.
\begin{Theorem}
\label{t4.3}
Assume that $\lambda>\lambda_{min}$ and that $\phi_1$, $\phi_2$ are two solutions 
of the ergodic problem \eqref{EP} associated to $\lambda$. If
\begin{equation}
\lim_{\vert x\vert\to+\infty}(\phi_1(x)-\phi_2(x))=0\, ,
\end{equation}
then $\phi_1=\phi_2$.
\end{Theorem}
This is once again in sharp contrast with the periodic case.

 \section{Entire Solutions of Hamilton-Jacobi Equations and Asymptotic Behavior}\label{ES}
 
 In this section we are interested in the solutions $v\in UC (\R^n\times \R)$\footnote{If $A\subset \R^m$, $UC(A)$ is the space of uniformly continuous functions on $A$.}  of Hamilton-Jacobi Equations set for all $t \in \R$, namely
 \begin{equation}\label{eqntt}v_t + F (x,Dv)=0\quad\hbox{in $\R^n\times \R $}\; .\end{equation}
 We are going to show that, under suitable conditions on $F$, $v$ is in fact independent of time, and is therefore solution of the stationary equation.
 Our key assumptions are
 \begin{itemize}
 \item[(H1)] There exists a viscosity subsolution $\phi \in UC(\R^n)$ of
 $F(x,D\phi )= 0$ in $\R^n$ such that $v-\phi$ is bounded.
 \item[(H2)] $F$ is bounded uniformly continuous in $\R^n \times B(0,R)$ for any $R>0$.
 \item[(H3)] There exists a continuous function
 $m:[0,+\infty )\rightarrow\R^+$
 such that $m(0^+)=0$ and, for all $x,y\in\R^n$ and $p\in\R^n$,
 $$|F(x,p)-F(y,p)|\le m(|x-y|(1+|p|))\; .$$
 \end{itemize}
 \noindent and,
 \begin{itemize}\item[(H4)] $\displaystyle{\left\{\begin{array}{l}\mbox{there exist $\eta >0$ and $\psi (\eta)>0$ such that, if $|F(x,p+q)|\ge\eta$ and}\\\mbox{$F(x,q)\le 0$ for some $x\in \R^n$, $p,q\in\R^n$, then, for all $\mu\in (0,1]$,}\\\\\hfill\mu F\Big( x,\mu^{-1}p+q\Big)\ge F(x,p+q)+\psi (\eta )(1-\mu).\hfill \end{array}\right.}$
 \end{itemize}
 It is worth noticing that if $F$ is $C^1$ in $p$, then (H4) reduces to
 \begin{itemize}\item[(H4)$'$] $\displaystyle{\left\{\begin{array}{l}\mbox{\hspace{1.0truein}$\displaystyle{F_p(x,p+q)\cdot p-F(x,p+q)\ge\psi (\eta)}$,}\\\\\mbox{for any $x\in \R^n$, $p,q\in \R^n$ such that $|F(x,p+q)|\ge \eta$ and $F(x,q)\le 0$.}\end{array}\right.}$
 \end{itemize}
 
 We show below that (H4) and (H4)' are satisfied if $F$ satisfies suitable (strong) convexity properties.
 
 The result about the solutions $u$ of (\ref{eqntt}) is the following
 \begin{Theorem}\label{ThmES}
 Assume that (H1)--(H4) hold.  Then any solution $v$ of (\ref{eqntt}) depends only on $x$, and is therefore a solution of $F=0$ in $\R^n$.
 \end{Theorem}
 Before commenting this result, we state its main consequence on the asymptotic behaviour of solutions of (\ref{HJ})-(\ref{idHJ}).
 \begin{Corollary}\label{ols}
 Assume that $H$ satisfies (\ref{coer}) and (H2), that $u_{0}$ is Lipschitz continuous in $\R^n$ and that there exists a solution $\phi$ of (\ref{EP}) for some $\lambda \geq \lambda_{min}$ such that $u_{0}-\phi$ is bounded in $\R^n$. If $H-\lambda$ satisfies (H4), then every function in the $\omega$-limit set of $u(\cdot,t)+\lambda t$ (in the sense of the local uniform convergence) is a solution of (\ref{EP}).
 \end{Corollary}
 
 This result may seem somewhat surprising, even in the bounded case, if we compare it to the counter-example of \cite{BS2} which shows that the (local uniform) convergence of $u(\cdot,t)+\lambda t$ as $t\to \infty$ may fail. It is worth pointing out anyway that Corollary~\ref{ols} does apply to this counter-example (since the nonlinearity satisfies strong convexity properties) and this demonstrates that Theorem~\ref{ThmES} is not sufficient to ensure such local uniform convergence as $t\to \infty$.  Again the problem we face here is the difference between local and global uniform convergence: under the assumptions of Corollary~\ref{ols}, if we have a sequence $(u(\cdot, t_n)+\lambda t_n)_n$ which converges uniformly in $\R^n$, we can conclude as in \cite{BS1} that $u(\cdot,t)+\lambda t$ converges as $t\to +\infty$, but this is wrong with only a local uniform convergence.
 
The assumptions (H4), (H4') are similar to the ones used in \cite{BS1}~: the only difference is that, on one hand, they concern here the whole set $\{|H|\geq \eta\}$ and not only the set $\{H \geq \eta\}$, and on the other hand, that it has to hold for $x$ in the whole space $\R^n$ while in \cite{BS1} several types of different behaviours can be mixed.
 
In \cite{BS1}, this assumption was a key condition to prove that, roughly speaking, a solution $v$ of such equation for $t\geq 0$, satisfies
$$ ||(v_{t})^-||_{\infty} \to 0 \quad\hbox{as  }t\to \infty\; ;$$
here this stronger formulation leads us to
$$ ||v_{t}||_{\infty} \to 0 \quad\hbox{as  }t\to \infty\; .$$
In order to understand why, we reproduce the formal argument provided in \cite{BS1} in the simpler case where $\phi=0$ - in fact, this formal argument is valid as soon as $\phi$ is $C^1$ and $D \phi$ is uniformly continuous in $\R^n$ - so that the transformation $\tilde v:=v-\phi$ can be done - and $F$ is smooth.

The Kruzhkov transform $w=-\exp(-v)$ provides a solution of
$$ w_{t} -w F(x,-\frac{Dw}{w}) = 0 \quad\hbox{in  } \R^n \times \R\; ,$$
and if we set $z=w_t$, it solves the linear equation
$$z_t+(F_p\cdot p-F)z+ F_p\cdot Dz=0\quad\hbox{in  } \R^n \times \R\; ,$$
where we have dropped the arguments of $F_p\cdot p-F$ and $F_p$ to simplify the notations. Next we consider $m(t)=\| z(\cdot,t)\|_\infty$; if $m(t) = z(x,t)$, we have $Dz(x,t)=0$ and the equation for $w$ implies that $\displaystyle z(x,t) = wF(x,-\frac{Dw}{w})$; therefore if $m(t) \geq \eta$, $F$ satisfies the same type of inequality and therefore (H4) says that $(F_p\cdot p-F) \geq \psi(\eta) >0$.  It then follows, that as long as $m(t) \geq \eta$
$$m'+ \psi(\eta) m = 0\; ,$$
 this implies that $m(t) \to 0 $ as $t\to \infty$.
 
The proof of Theorem~\ref{ThmES} will make this formal proof more precise~; since it is very similar to the proof in the Appendix of \cite{BS1}, we will just sketch it, pointing out the main adjustments. Now we check assumption (H4) ; a typical case we have in mind is the case when, on one hand, we consider Lipschitz  continuous solutions, and, on an other hand, $F$ is $C^2$ in $p$ for any $x$ and satisfies, for some $\beta >0$
$$ F_{pp}(x,p) \geq \beta \hbox{Id} \quad\hbox{in  } \R^n \times \R^n\; .
$$
By the convexity of $F$, one has $$ \mu F(x, \mu^{-1}p+q)-F(x,p+q) \geq -(1-\mu) F(x,q) + \frac{\beta}{2}\mu(1-\mu)|p|^2\; .$$
 Since
we consider Lipschitz continuous solutions we can assume that $|p|, |q|\leq K$ for some large enough constant $K$ and, thanks to assumption
(H2), there exists a modulus of continuity $m$ for $F$ on $\R^n \times B(0,K)$.  Now, assume that $|F(x,p+q)|\geq \eta$ and $F(x,q) \leq 0$; if $F(x,q) \geq -\eta/2$, we have at the same time $|F(x,p+q) - F(x,q)| \leq m(|p|)$ and $|F(x,p+q) - F(x,q)| \geq \eta/2$ and therefore $|p| \geq \chi(\eta) >0$. And (H4) is satisfied because of the term $\di\frac{\beta}{2}\mu(1-\mu)|p|^2$, the other one being positive.
If on the contrary, $F(x,q) \leq -\eta/2$, then the term $-(1-\mu) F(x,q)$ provides the positive sign. 

These computations shows that (H4) is related to the strict convexity of $F$. Strict convexity is not optimal, but the counterexample 
$$ v_t + |v_x +\alpha| - |\alpha | = 0 \quad\hbox{in  } \R \times \R\; $$
analyzed in \cite{BS1}, shows that $H$ must not be too far from being strictly convex. 

\noindent{\bf Proof of Corollary~\ref{ols}}. Since $u_0 - \phi$ is bounded and since $\phi(x) - \lambda t$ is a solution of (\ref{HJ}), the comparison principle for viscosity solutions yields
$$ ||u - \phi(x) + \lambda t||_\infty \leq || u_0 - \phi ||_\infty\; .$$
We set $\tilde u (x,t) = u(x,t) + \lambda t$. By the above inequality, $\tilde u - \phi$ is uniformly bounded and $\tilde u$ is a solution of
$$ \tilde u_t + F(x,D\tilde u) = 0 \quad\hbox{in  } \R^n\times (0, +\infty)\; ,$$
where $F(x,p) = H(x,p) -\lambda$ in $\R^n \times \R^n$.

If $w$ is in the $\omega$-limit set of $\tilde u$, then there exists a sequence $(t_p)_p$ converging to $+\infty$ such that $\tilde u (\cdot, t_p)$ converges locally uniformly to $w$.
We set $v_p (x,t) = \tilde u(x,t_p+t)$. The function $v_p$ is a viscosity solution of
$$ (v_p)_t + F(x,Dv_p) = 0 \quad\hbox{in  } \R^n\times (-t_p, +\infty)\; ,$$
and extracting if necessary a subsequence (since $v_p$ is uniformly bounded and Lipschitz continuous), we may assume that $v_p$ converges locally uniformly to Lipschitz continuous function $v$, defined on $\R^n \times \R$, such that $v-\phi$ is bounded. Moreover by stability result for viscosity solutions, $v$ solves (\ref{eqntt}). 

Since $\phi$ is a solution of $F=0$, we can use Theorem~ \ref{ThmES} and deduce that $v$ is independent of $t$ and is a solution of $F=0$. But, since $v(x,t)=v(x,0) = w(x)$, we have proved that $w$ is a solution of $H=\lambda$. \hfill $\bullet$

\medskip

\noindent {\bf Proof of Theorem~ \ref{ThmES} (sketch).}  Changing $\phi$ in $\phi-C$ for some constant $C>0$ large enough, we can assume that $v-\phi \geq 1$ in $\R^n \times \R $. 
Then for $\eta >0$, we introduce the functions$$\mu^+_\eta (t)=\max_{x\in\R^n, s\ge t} \Big[\frac{v(x,s)-\phi (x)+2\eta(s-t)}{v(x,t)-\phi (x)}\Big]\; , $$
$$
\mu^-_\eta (t)=\min_{x\in\R^n, s\ge t} \Big[\frac{v(x,s)-\phi (x)-2\eta(s-t)}{v(x,t)-\phi (x)}\Big]\; .$$
The functions $\mu^+_\eta$, $\mu^-_\eta : \R \to\R$ are
Lipschitz continuous, and $\mu^+_\eta \geq 1$, $\mu^-_\eta \leq 1$ in $\R$. The key point of the proof is to show that $\mu^+_\eta$ is a subsolution of the variational inequality
$$ \min\left(\mu' (t)+k\psi (\eta )(\mu (t)-1), \mu(t)-1\right) =  0\quad\mbox{in $\R$}\; ,$$
and that $\mu^-_\eta$ is a supersolution of the variational inequality
$$\max\left( \mu' (t)+k\psi (\eta )(\mu (t)-1), \mu(t)-1\right) =  0\quad\mbox{in $\R$}\; ,$$
for some constant $k>0$. The proof of this fact for $\mu^-_\eta$ is given in the Apendix of \cite{BS1} and for $\mu^+_\eta$, the proof is almost the same with few minor changes.

We deduce from these properties and the uniqueness results for these variational inequalities, that $\mu^+_\eta(t) ,  \mu^-_\eta(t) \equiv 1$; indeed by choosing $T>0$ large
$$ |\mu^+_\eta(t) -1|, |\mu^+_\eta(t) -1| \leq \tilde C \exp (-k\psi (\eta )T)\; ,$$
where $\tilde C= \max(| |\mu^+_\eta(t) -1||_\infty, ||\mu^+_\eta(t) -1||_\infty) $. And letting $T \to \infty$ provides the results.

This equality being valid for any $\eta$, we deduce  that $v$ is independent of time. \hfill$\bullet$

\end{document}